\documentclass{article}

\usepackage{makeidx}

\usepackage{array}

\usepackage{natbib}
\usepackage{graphicx}
\usepackage{booktabs}

\usepackage{url}

\usepackage{amsmath,amsfonts}

 \newcommand{\Expect}[1]{\operatorname{\mathbb{E}}\left[#1\right]}
 \newcommand{\Indicator}[1]{\operatorname{\mathbb{I}}\left[#1\right]}

 \newcommand{\Prob}[1]{\operatorname{\mathbb{P}}\left[#1\right]}

 \newcommand{\Var}[1]{\operatorname{Var}\left[#1\right]}

 \DeclareMathOperator*{\argmin}{arg\,min}
 
 \renewcommand{\d}{\operatorname{d}}

 \newcommand{\dist}{\operatorname{dist}}
 \newcommand{\half}{\tfrac{1}{2}}
 
 \newcommand{\re}{\operatorname{e}}
 \newcommand{\botrule}{\bottomrule}
\begin{document}

%
%
%
%
%
%
%
%
\title{Barycentres and Hurricane Trajectories}%
\author{Wilfrid S.~Kendall, University of Warwick\footnote{Work supported by EPSRC grant EP/K031066/1}}
\date{\text{~}}
\maketitle

\begin{abstract}\noindent
The use of barycentres in data analysis is illustrated, using as example a dataset of hurricane trajectories.
\end{abstract}


\noindent
\textbf{Keywords and phrases:}\\
\textsc{Barycentre; Cosine barycentre; One-dimensional Potts model; \\
Riemannian barycentre; Run statistic.}

\section{Introduction}\label{sec:WSKintro}
This paper is principally motivated by intellectual curiosity. 
After work with Huiling Le \citep{KendallLe-2011} on laws of large numbers and central limit theorems for empirical Riemannian barycentres, 
it seemed natural to investigate the use of Riemannian barycentres in data analysis.
This topic relates to contemporary interest in the statistical analysis of data comprised of intrinsically geometric objects, which can be viewed as part of the subject sometimes known as `object-oriented data analysis'.
Indeed the statistical use of barycentres has already been pioneered, for example, in \citet{FournelReynaudBrammerSimmonsGinestet-2012} and \citet{GinestetSimmonsKolaczyk-2012a} (also see early work by \citealp{Ziezold-1994});
the purpose of the present paper is to use a specific application to explore their use in analyzing trajectories with strong geometric content.
In the following we use Riemannian barycentre theory to produce a simple non-parametric analysis of the extent to which consecutive North Atlantic hurricanes might have similar behaviour.
Note that considerably more sophisticated methods of curve-fitting on manifolds could have been used here (see for example the use of smoothing splines described in \citealp{SuDrydenKlassenLeSrivastava-2011}): 
the barycentre approach is relatively simplistic, but nonetheless may be useful.

Writing this paper affords the opportunity of expressing sincere homage to Kanti Mardia for his seminal leadership in the application of geometry to statistics.
I owe him thanks for kindness and encouragement stretching right back to 1978, when as a callow research student I was invited by Kanti to make a most stimulating research visit to the
University of Leeds Statistics group. 
Moreover, the present work originated in preparation for a talk I gave at one of the famous LASR workshops initiated and nurtured by Kanti at Leeds (specifically LASR 2011).
I hope that Kanti finds pleasure in reading this brief account.

The paper commences (Section \ref{sec:WSKbarycentres}) with a speedy review of relevant aspects of Riemannian barycentres, with special attention paid to the simple but fundamental case of sphere-valued data.
Section \ref{sec:WSKhurricanes} then reviews \texttt{HURDAT2}, 
a remarkable publicly available dataset composed of hurricane trajectories (tropical cyclones in the North Atlantic) and some associated data
concerning wind speeds and atmospheric pressures. 
This is the test dataset: attention will be confined to the hurricanes viewed as trajectories on the terrestrial sphere.
We then describe (Section \ref{sec:WSKkmeans}) how barycentre theory interacts with non-parametric statistics \emph{via} \(k\)-means clustering, and discuss preliminary results for the test dataset (Section \ref{sec:WSKresults}).
The concluding Section \ref{sec:WSKconc} reviews the results and considers some possible next steps.

\section{Barycentres}\label{sec:WSKbarycentres}
\citet{Frechet-1948} introduced barycentres in metric spaces as minimizers of `energy functionals' \(x\mapsto\Expect{\dist^2(X,x)}\) for random variables \(X\) taking values in metric spaces.
\citet{Kendall-2013c} presents a recent review of some subsequent theory; there is a strong link with convexity \emph{via} `convex geometry'\index{convex geometry} (\citealp{Kendall-1990d}, see also \citealp{Afsari-2010}).
Our interest is focussed on the theory of Riemannian barycentres\index{barycentre!Riemannian barycentre} for random variables taking values in the \(2\)-sphere \(S^2\): indeed this can be viewed as normative for Riemannian barycentres \citep{Kendall-1991a}.
In particular, sphere-valued random variables can be guaranteed to possess unique barycentres when their distributions are concentrated in closed subsets of open hemispheres (which is to say, when the random variables are
confined to `small hemispheres').
Considerable work has been devoted to establishing laws of large numbers and central limit theorems for empirical barycentres \citep{PatrangenaruBhattacharya-2003,BhattacharyaPatrangenaru-2005,BhattacharyaBhattacharya-2008};
this has even opened up a new multivariate perspective on the classical Feller-Lindeberg central limit theory \citep{KendallLe-2011}.
An interesting non-Riemannian case is discussed in the pioneering work of \citet{HotzHuckemannLeMarronMattinglyMillerNolenOwenPatrangenaruSkwerer-2012};
see also \citet{BardenLeOwen-2013}.
In the present paper our interest centres on more data-analytic concerns, based on barycentres of measurable random maps \(\Phi:[0,T]\to S^2\) from a time interval to \(S^2\). 
Convex geometry, hence uniqueness of barycentres, is maintained if for each time \(t\) the random variable \(\Phi(t)\) is supported in a (possibly time-varying) small hemisphere.

There are a number of studies of iterative algorithms for computing Riemannian barycentres (for example, \citealp{Le-2004,ArnaudonDombryPhanYang-2010}).
We shall finesse such considerations by approximating Riemannian barycentres on the sphere by \emph{cosine barycentres}\index{barycentre!cosine barycentre}, given by
projecting the conventional expectation onto the sphere,
\[
 \argmin_{x\in S^2} \Expect{1 - \cos\dist(X,x)} \quad=\quad \frac{\Expect{X}}{\|\Expect{X}\|}\,.
\]
This is the `mean direction' in the terminology of directional statistics. We choose to use the term ``cosine barycentre'', to emphasize that it
minimizes what one might term the \emph{cosine-energy} (related to chordal distance) \(1 - \cos\dist(x,y) \approx \half \dist(x,y)^2\). 
In particular, it provides a feasible and explicit approximation to the Riemannian barycentre if the dispersion of \(X\) on \(S^2\) is not too large.
Note however that its explicit and constructive definition provides no easy panacea for questions of uniqueness: evidently the construction only works when \(\Expect{X}\neq0\), 
and indeed if the support of the data cannot be contained in a closed subset of an open hemisphere then it is possible for the barycentre to be ill-defined. 
(Consider, for example, the problem of finding the cosine-barycentre for a probability distribution spread out uniformly over a fixed great circle. 
Then all points on the sphere minimize the cosine-energy.)

\begin{figure}
 \includegraphics[width=15cm]{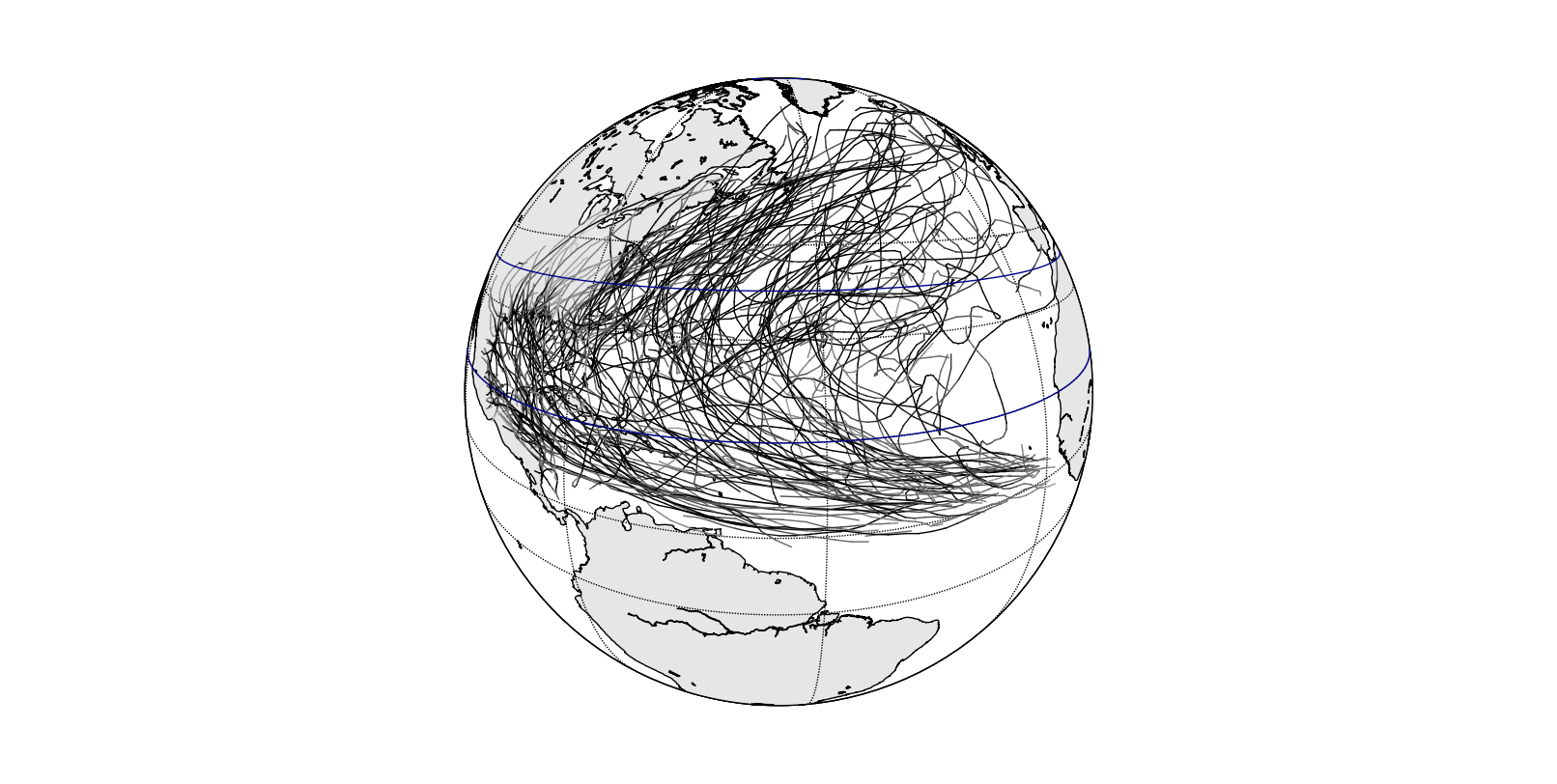}
 \caption{Trails of the 233 hurricanes recorded over the period 2000-2012 in the \texttt{HURDAT2} dataset
 (darker trajectories indicate high maximum sustained wind speed).
 The viewpoint of this and similar following images is placed \(5000\) km above the centre of the image,
 which is therefore distorted at normal viewing distance for all but the most extremely short-sighted. 
 We will consider the longer period 1950-2012, and will restrict attention to hurricanes crossing \(20^\circ N\) and \(35^\circ N\) (drawn as continuous lines in figure), 
 and register hurricanes on their first crossing of \(35^\circ N\).
 }
 \label{fig:WSKfig01}
\end{figure}

Cosine-barycentres allow us to choose representative barycentre trajectories\index{trajectory!barycentre trajectory} for a collection of hurricane trajectories\index{trajectory!hurricane trajectory},
defined as solving the minimization problem
\begin{equation}\label{eq:barycentre}
 \argmin_{F:[0,T]\to S^2} \frac{1}{T}\int_0^T \Expect{1 - \cos\dist(F(t), \Phi(t))} \d t\,.
\end{equation}
Here the expectation is actually the empirical sample average over the collection of hurricane trajectories,
so that \(\Phi\) is viewed as drawn uniformly at random from this collection. Note that the minimization
in \eqref{eq:barycentre} can be carried out separately for each time \(t\), since there is no continuity requirement placed on \(F\).
We choose to avoid considerations of continuity or of smoothness of trajectories; 
close inspection of the hurricane trajectories in Figure \ref{fig:WSKfig01} suggest that smoothness, at least, is perhaps not a paramount consideration.
However there are some practical issues that need to be faced. 
Our hurricane trajectories are not \emph{a priori} registered to comparable starting and / or finishing times; in fact typically their durations do not overlap.
We restrict attention to trajectories which make upcrossings on latitudes of \(20^\circ N\) and \(35^\circ N\): we register times to agree at
the first upcrossing of latitude \(35^\circ N\).
(The following analysis is sensitive to these choices: lower latitudes do not produce a clear statistical signal.)

A further issue is the need to average over hurricanes which, even when time-registered, start at different negative times and finish at different positive times.
A possible solution is suggested by the observation that all hurricanes under consideration are confined to the Northern hemisphere.
This suggests the following idea: if \(\Phi_\omega(t)\) is not defined for some sample point \(\omega\in\Omega\) then replace \(1 - \cos\dist(F(t), \Phi_\omega(t))\) in \eqref{eq:barycentre}
by the average when \(\Phi(t)\) is uniformly distributed over the equator
(this decision could be justified by arguments of maximum entropy type).
However,  this modification does not lead to good results in our current application.
Instead we solve the issue by restricting attention to the largest time interval over which our collection of time-registered hurricanes all have defined locations.
This cropping procedure has the disadvantage of restricting the analysis to behaviour of the trajectories near the specified upcrossing used for registration.

Finally, we mention two more possible refinements, both of which would be  computationally demanding and which we will not adopt.
Firstly, one could attempt to solve the extended minimization problem allowing time-shifts of individual trajectories within the minimization problem \eqref{eq:barycentre}.
Secondly, one could replace the time integral in \eqref{eq:barycentre} by an integral over arc-length, or perhaps an integral over upcrossings of latitudes 
(though in this case one would be deliberately introducing discontinuity in cases when hurricane trajectories decrease as well as increase in latitude).
Indeed one could envisage a whole variety of possible nonlinear warpings of time.
We defer to another occasion the consideration of these refinements, as well as of assessment of the effect of approximating Riemannian barycentres by cosine-barycentres.

Our application needs to find a number of representative barycentre trajectories instead of just one.
The natural remedy is to apply the \(k\)-means algorithm\index{k-means algorithm@\(k\)-means algorithm}, adapted to use cosine-barycentres.
Specifically, we aim to use barycentre \(k\)-means to cluster the chosen set of hurricane trajectories, 
so that we can study temporal association between cluster labels defined by the barycentre trajectory to which each hurricane trajectory is attached.
We use Lloyd's algorithm \citep{Lloyd-1982} for \(k\)-means\index{Lloyd's algorithm|see {\(k\)-means algorithm}}: beginning with a random
initial set of \(k\) trajectories serving as cluster centroid trajectories,
the algorithm alternates between associating each trajectory to the closest cluster centroid trajectory (measured by cosine-distance), 
and then replacing each cluster centroid trajectory by the computed barycentre trajectory for the cluster.
The algorithm has to be run repeatedly in order to find a good clustering; we choose to use \(10\) repetitions. 
Typically a single run of the algorithm will not produce an optimal set of cluster centroid trajectories (here, minimizing within-cluster sum of cosine-distances); 
indeed the task of producing such an optimal set is typically NP-complete.

Faster algorithms \emph{do} exist for the one-dimensional problem \citep{WangSong-2011}, 
but it is an open problem to extend these to `nearly one-dimensional' structure as exhibited by the set of hurricane trajectories.
In this study we use the \(k\)-means algorithm, setting \(k=20\), to group hurricanes into \(20\) groups linked to \(20\) barycentre trajectories. 
The groups are ordered from west to east according to where the barycentre trajectories first cross
latitude \(35^\circ N\).

\section{Hurricanes}\label{sec:WSKhurricanes}
As an illustrative application, we consider the remarkable and freely available \texttt{HURDAT2} dataset\index{HURDAT2 dataset@\texttt{HURDAT2} dataset} 
concerning hurricanes (tropical cyclones) of the North Atlantic ocean, 
a collection of 1740 hurricane trajectories in the North Atlantic recorded by various means over 161 years from 1851 to 2012
(see Figure \ref{fig:WSKfig01} for a display of recent hurricanes; the dataset is discussed in \citealp{LandseaFranklin-2013,McAdieLandseaNeumannDavidBlakeHammer-1978}). 
The dataset is available at
\begin{center}
 \url{www.aoml.noaa.gov/hrd/hurdat/Data_Storm.html}.
\end{center}
Our interest centres on whether there is any evidence for temporal association; is there a tendency for successive hurricane trajectories to be close?
\citet{MacManus-2011} used geometric methods to investigate similar issues (area of overlapping curvilinear strips based on paths, and non-parametric measures of association).
Here, we intend to use this question to illustrate application of the notion of barycentres based on hurricane paths considered as \(S^2\)-valued trajectories.
Evidently the extent of these paths, ranging over wide expanses of the North Atlantic, means that their underlying geometric nature should be taken seriously.

We emphasize that this investigation is of purely methodological interest, aimed at illustrating the use of barycentres in data analysis. 
Addressing the question of temporal association properly would require serious attempts to relate \texttt{HURDAT2} to other datasets and sources of information.

\begin{figure}
 \includegraphics[width=15cm]{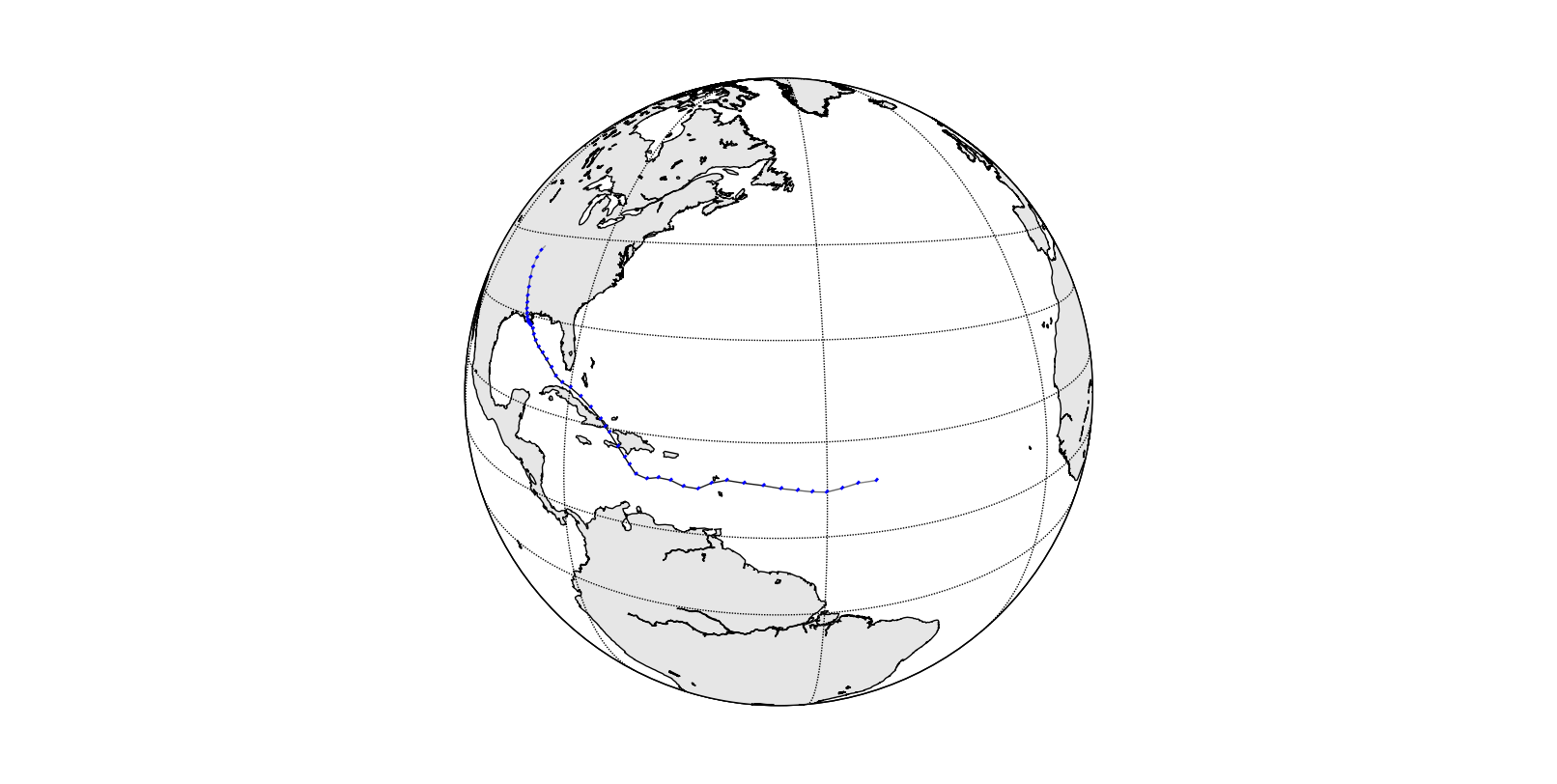}
 \caption{Sampling points (measured every 6 hours) of the Best Track of hurricane Isaac in 2012.
Typical separation between measurement points is about 100km.}
 \label{fig:WSKfig02}
\end{figure}

\begin{figure}
 \includegraphics[width=10cm]{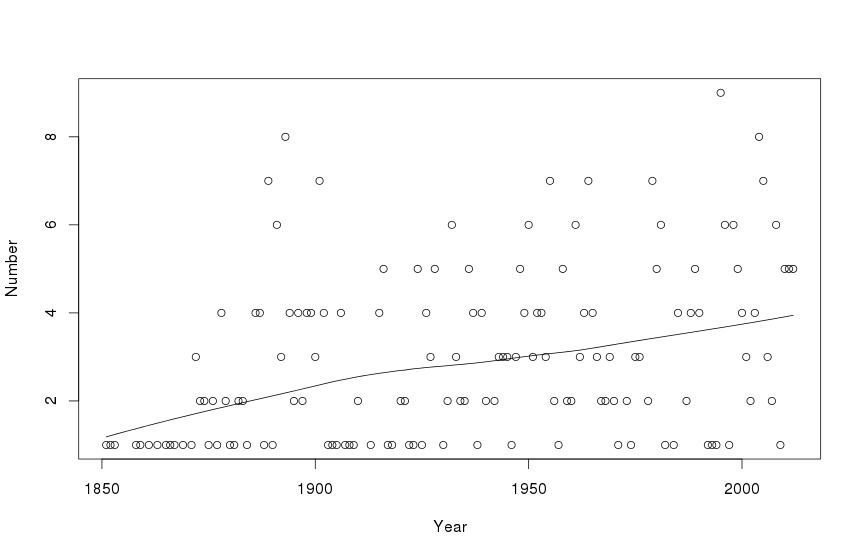}
 \caption{Numbers of hurricanes per year in the \texttt{HURDAT2} dataset over 1851-2012 which cross latitudes \(20^\circ N\) and \(35^\circ N\), together with fitted \texttt{lowess} curve.}
 \label{fig:WSKfig0n}
\end{figure}

\begin{figure}
 \includegraphics[width=10cm]{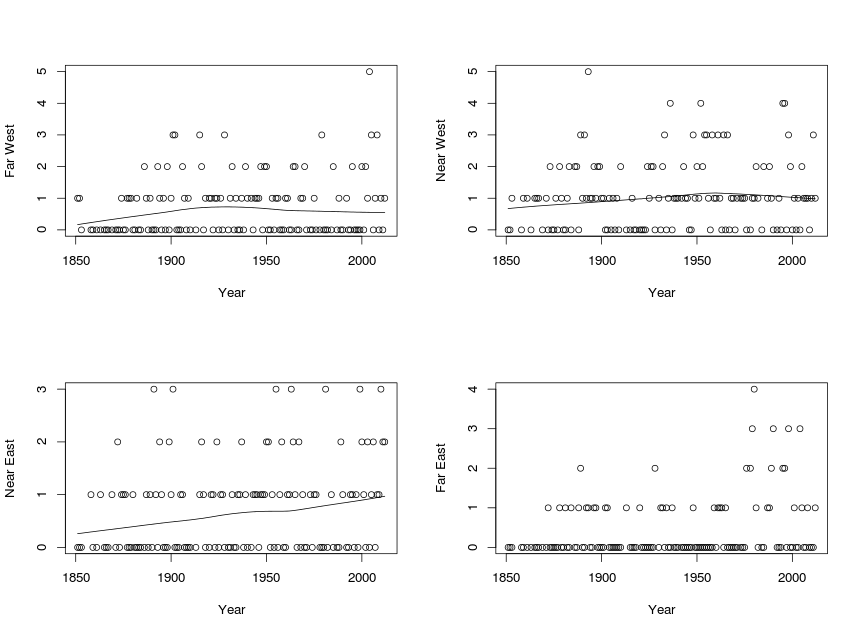}
 \caption{Numbers of hurricanes per year in the \texttt{HURDAT2} dataset over 1851-2012 which cross latitudes \(20^\circ N\) and \(35^\circ N\), grouped according to whether they 
 first cross latitude \(35^\circ N\) at far-west, near-west, near-east, or far-east locations as indicated using \(k\)-means clustering (with \(k=4)\)
 together with fitted \texttt{lowess} curves (except in the far-east case, for which the large majority of years record no hurricanes).}
 \label{fig:WSKfiga4}
\end{figure}

The following remarks are taken from  \citet{McAdieLandseaNeumannDavidBlakeHammer-1978}
(further useful background is also given in a striking statistical survey of data 
acquired from three years of USAF flight missions flown into Northwest Pacific tropical cyclones\index{tropical cyclone|see {hurricane}} which is reported by \citealp{WeatherfordGray-1988a,WeatherfordGray-1988b}).
Each hurricane trajectory in the \texttt{HURDAT2} dataset is a `Best Track',\index{hurricane!`Best Track'} defined using best estimates of the location of the hurricane centre\index{hurricane}, 
reconciling measurements obtained by various means, and taken at 6-hourly intervals with small-scale smoothing applied. 
In essence each hurricane trajectory is represented by a timed sequence of latitude / longitude pairs, measured every quarter of a day in degrees of latitude and longitude to an accuracy of 1 decimal place. 
Recall that a degree of latitude or longitude at the equator represents separation slightly in excess of 110 km, so stated location precision is of order \(\pm5\) km.
For the North Atlantic basin, \citet{McAdieLandseaNeumannDavidBlakeHammer-1978} reports that diameters of hurricane eyes lie in the range 15-50 km, so this is an entirely adequate level of precision. 
Typically hurricane eyes move at speeds of order 5 m/s, so successive six-hourly measurements are separated by order of 100 km (Figure \ref{fig:WSKfig02}).
Measurements of maximum sustained surface wind-speed and (subsequent to 1979) central surface barometric pressure are also recorded at various distances from hurricane centres; 
however we focus on geographic location alone.

Surveillance methods for data capture have of course varied over the period of the dataset. Early observations were acquired from sailing ship logs; in due course these were supplemented by radio reports, 
then by aircraft and radar observations, and finally by satellite observation and other systems such as dropsondes.
Whether because of increased observational capacity, or whether because of secular change (long-term non-periodic variation) in global weather conditions,
the number of recorded hurricanes over the period 1851-2012 is clearly increasing with time (see the \texttt{lowess} curve for 
annual counts of hurricanes crossing \(20^\circ N\) and \(35^\circ N\), given in Figure \ref{fig:WSKfig0n}): 
for example, later years seem more likely to experience 6 or more such hurricanes.
Further evidence of secular change is obtained by categorizing these hurricanes using the \(k\)-means algorithm, with \(k=4\), applied to those hurricanes crossing latitudes \(20^\circ N\) and \(35^\circ N\).
Figure \ref{fig:WSKfiga4} indicates that the east-west distribution of recorded hurricanes also varies over this time period; it appears that more easterly hurricanes are recorded in later years.

It is evident from this discussion that, for the purposes of study of temporal association, the \texttt{HURDAT2} dataset is best considered as a 
temporally-ordered sequence of short time-series, one short time-series per year over 1851-2012;
moreover the statistics of these short time series should be expected to be different in later as opposed to earlier years.
We shall focus below on the period 1950-2012, but should bear in mind the trends illustrated in Figure \ref{fig:WSKfiga4}.

\section{Using $k$-means and non-parametric statistics}\label{sec:WSKkmeans}
We seek to employ barycentre techniques to investigate temporal association between successive hurricanes.

Regardless of the \emph{ad hoc} aspects of clustering using the \(k\)-means algorithm (specifically, a potential dependence 
of actual implementation outcomes on initial conditions),
the crucial point is that the clustering takes no account of time order, whether by year or by time within year.
Consequently this clustering method can be used to detect temporal association using non-parametric statistical permutation tests\index{non-parametric statistical permutation test}. 

Note that we do not consider the issue of estimation of the number of clusters \(k\), as we are interested in the output of the \(k\)-means algorithm solely
as an intermediate device to facilitate the detection of temporal association.

The simplest choice is to compute the statistic \(T\)
counting the number (summed over all years)
of pairs of hurricanes, one immediately succeeding the other in a given year, such that both hurricanes are categorized as belonging to the same \(k\)-mean cluster.
Note that this is effectively a multiple category variant of the runs test\index{Wald-Wolfowitz runs test for randomness} for randomness, as discussed by \citet{WaldWolfowitz-1940}, 
since each hurricane either initiates a run, or completes a successive pair of hurricanes from the same cluster.
(A discussion of the multiple category variant is given by \citealp{Mood-1940}, who also records an informative early history of the runs test.)

It is enlightening to impose a narrative based on an informal statistical model of temporal variation,
as this allows a somewhat more statistically principled approach and suggests a useful generalization.
Independently for each year \(y\), let \(X_{y,i}\) (for \(i=1, 2, \ldots, h_y\)) record the \(k\)-mean cluster containing the \(i^\text{th}\) hurricane (measured in time-order) of the \(h_y\) hurricanes observed that year. 
Then the year \(y\) in question contributes the following summand \(T_y\) to the total non-parametric statistic \(T=\sum_y T_y\):
\[
 T_y(x_1, \ldots, x_{h_y})\quad=\quad\sum_{i=2}^{h_y} \Indicator{X_{y,i}=X_{y,i-1}}\,.
\]
This suggests that we model the sequence of \(k\)-mean labels in a given year by a one-dimensional Potts model\index{Potts model} \citep[\S1.3]{Grimmett-2006}: conditional on \(h_y\) the total number of hurricanes in that year,
the probability of observing the sequence \(x_1, x_2, \ldots, x_{h_y}\) is
\begin{multline}\label{eq:potts}
\Prob{X_{y,1}=x_1, \ldots, X_{y,h_y}=x_{h_y}| h_y}\quad=\quad\\
 \frac{1}{Z(\theta, h_y)}\exp(\theta\;T_y)\quad=\quad
 \frac{1}{Z(\theta, h_y)} \exp\left(\theta\;\sum_{i=2}^{h_y} \Indicator{x_i=x_{i-1}}\right)\,.
\end{multline}
Here \(\theta\geq0\) is the parameter relating to strength of association, and the partition function can be computed explicitly in this simple one-dimensional case: \(Z(\theta, h)=k (\re^\theta + k - 1)^{h-1}\).
(The \(k=2\) case corresponds to a one-dimensional Ising model\index{Ising model}, and can be traced as far back as Ernst Ising's 1924 thesis, published in part in Zeitschrift f\"ur Physik in 1925.)

Treating each year \(y\) as independent, we choose to condition not only on the total number \(h_y\) of hurricanes in that year, 
but also on the total numbers \(R_{y,j}=\#\{i:X_{y,i}=j\}\) of hurricanes in the year \(y\) categorized as belonging to \(k\)-mean cluster \(j\).
The action of conditioning on the \(R_{y,j}\) discards some information about \(\theta\), since large positive values of \(\theta\) would promote dominance by a single cluster in each year \(y\).
However  the empirical evidence of secular trends supplied by Figure \ref{fig:WSKfiga4}
suggests the need for a more realistic model for the measurements \(R_{y,j}\), allowing for their distributions not being exchangeable over the label \(j\);
conditioning on the \(R_{y,j}\) allows us to evade this difficulty.
The score statistic for the resulting conditioned model at \(\theta=0\) is then the year's contribution \(T_y\) to the non-parametric statistic \(T\).
Thus consideration of \(T\) amounts to performing a conditional Neyman-Pearson hypothesis test of a null  hypothesis \(\mathcal{H}_0: \theta=0\), 
against a one-sided compound hypothesis \(\mathcal{H}_1: \theta>0\).

It is of course possible to develop this theme further, using the evaluation of \(Z(\theta, h_y)\) for the one-dimensional Potts model. Thus:
  \begin{enumerate}
   \item Inference could be improved to take explicit account of the information provided by the pattern of values of \(R_{j,y}\),
   for example, by imposing an external field on the Potts model \eqref{eq:potts} to obtain
\begin{multline*}
\Prob{X_{y,1}=x_1, \ldots, X_{y,h_y}=x_{h_y} | h_y}\quad\propto\quad   \\
\exp\left(\theta\;\sum_{i=2}^{h_y}\Indicator{x_i=x_{i-1}} + \sum_{i=1}^{h_y}\sum_{j=1}^k \psi_j \Indicator{x_i=j}\right)\,;
\end{multline*}
   \item Or one could attempt maximum likelihood estimation of \(\theta\), or even of different \(\theta_{j}\) pertaining to different clusters \(j\) using a refined probability mass distribution
\begin{multline*}
\Prob{X_{y,1}=x_1, \ldots, X_{y,h_y}=x_{h_y} \;|\; h_y, R_{y,1}, \ldots, R_{y,k}}\quad\propto\quad\\
\exp\left(\sum_{i=2}^{h_y} \theta_{x_i} \Indicator{x_i=x_{i-1}}\right)
 \,.
\end{multline*}
  \end{enumerate}
However we avoid pursuing either of these leads here, not only because of the resulting increase in model complexity but also because this would commit us 
in excessive detail to a parametric model which is motivated largely by heuristic considerations. Moreover the model is dependent on the output of the \(k\)-means Lloyds algorithm, itself potentially a random phenomenon, insofar as the actual outcome of the algorithm can depend on essentially random selection of initial conditions (mitigated by using repeated independent runs, and taking the best resulting clustering).
Finally, the model, as expressed here, uses (na{\"\i}ve) free boundary conditions (it treats initial and final hurricanes of each year in much the same way as all the others), 
which is a further reason not to take it too seriously.

To evaluate the significance of the score statistic \(T\), we compute the conditional mean and variance of each \(T_y\), add up the \(T_y\) over the years \(y\) under consideration,
and refer the sum \(T\) to a normal distribution with matching mean and variance (this relies implicitly on a central limit theorem approximation of Lyapunov type),
or use a simulation test based on random permutations within each year.
Means and variances can be computed by straightforward combinatorial methods: suppose that in year \(y\) there are present \(R_{y,j}=r_j\) hurricanes belonging to cluster \(j\), for \(j=1, \ldots, k\).
For convenience we set \(m_2=\sum_{j=1}^k r_j(r_j-1)\) and \(m_3=\sum_{j=1}^k r_j(r_j-1)(r_j-2)\), and find
\begin{align}
 \Expect{T_y | R_{y,1}=r_1, \ldots, R_{y,k}=r_k} \quad&=\quad \frac{m_2}{h_y}       \,,\label{eq:mean}\\
 \Var{T_y | R_{y,1}=r_1, \ldots, R_{y,k}=r_k}    \quad&=\quad \frac{m_2^2}{h_y^2(h_y-1)} + \frac{(h_y-3) m_2}{h_y(h_y-1)} - \frac{2 m_3}{h_y(h_y-1)}\,.\label{eq:variance}
\end{align}
Note that differentiation of the partition function \(Z(\theta, h_y)\) would produce means and variances \emph{not} conditioned on the pattern of \(R_{y,j}\) values,
which would not suit our purpose.

It is useful to modify the one-dimensional Potts model \eqref{eq:potts} to allow for a geometric decay in strength of association, with decay rate \(\beta\in(0,1)\) and \(\beta\ll1\):
\begin{multline}\label{eq:potts-spread}
\Prob{X_{y,1}=x_1, \ldots, X_{y,h_y}=x_{h_y}| h_y}\quad=\quad\\
 \frac{1}{Z(\theta, h_y, \beta)}\exp(\theta\;T_{y,\beta})\quad=\quad
 \frac{1}{Z(\theta, h_y, \beta)} \exp\left(\theta\;\sum_{\ell=1}^{h_y-1} \beta^{\ell-1} \sum_{i=\ell+1}^{h_y} \Indicator{x_i=x_{i-\ell}}\right)\,.
\end{multline}
The corresponding conditional score statistic (given the pattern of values \(R_{y,j}\)) is 
\begin{equation}\label{eqn:potts-spread-stat}
 T_{y,\beta}\quad=\quad 
  \sum_{\ell=1}^{h_y-1} \beta^{\ell-1} 
  \sum_{i=1+\ell}^{h_y}\Indicator{X_{y,i}=X_{y,i-\ell}}
 \quad=\quad
 \sum_{\ell=1}^{h_y-1} \beta^{\ell-1} 
\#\left\{\text{agreements at lag }\ell\right\}\,,
\end{equation}
and this provides some capacity to allow for detection of temporal association between trajectories which are not immediately consecutive in time.
The computations of mean and variance at \eqref{eq:mean} and \eqref{eq:variance} can be generalized to cover this case, 
though the formulae are too unwieldy to report here.
In any case, in the sequel we shall evaluate test statistic scores using a simulation test based on random permutations within each year.

\section{Results}\label{sec:WSKresults}
\begin{figure}
 \includegraphics[width=15cm]{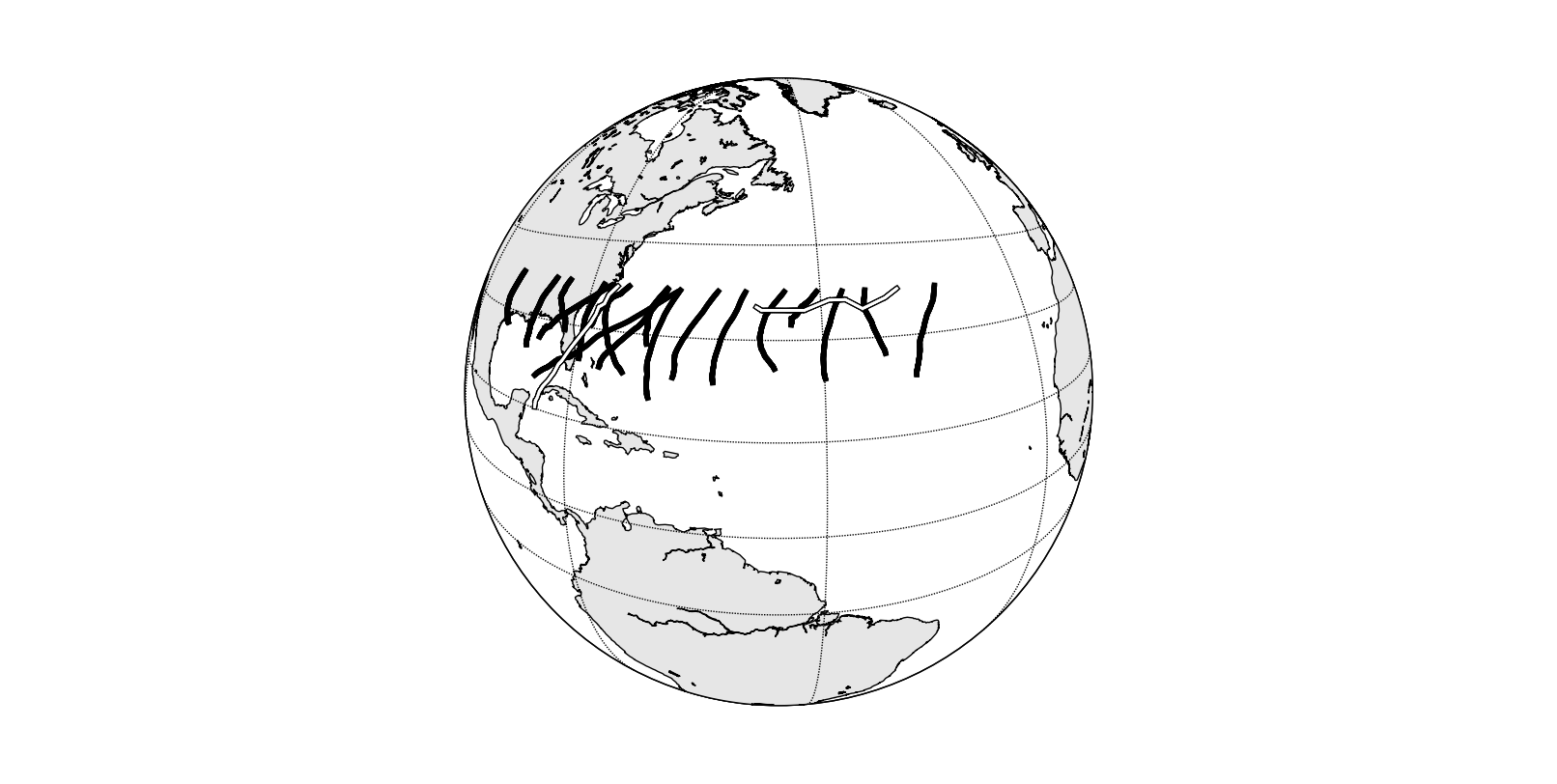}
 \caption{Plot of \(20\) barycentre trajectories arising from the \(k\)-means algorithm with \(k=20\), applied to the 1950-2012 dataset of hurricanes crossing \(20^\circ N\) and \(35^\circ N\).
 Barycentre trajectories are denoted by thick paths (the two outline paths correspond to clusters of just one or two hurricanes).
 The apparently anomalous trajectory running nearly horizontally near the more eastern end of the collection arises from a single rather long trajectory, whose initial behaviour 
 (including its up-crossing of \(20^\circ N\)) is cropped as part of the process of cropping all hurricane trajectories to the same maximal time-interval.
 }
 \label{fig:WSKfig04}
\end{figure}
We consider a specific example, namely 179 North Atlantic hurricanes occurring over the time-range 1950-2012, and crossing latitudes \(20^\circ N\) and \(35^\circ N\),
and registered at first upcrossing of \(35^\circ N\). 
It is candidly admitted that other crossing latitude choices lead to results which are not statistically significant, so assessment
of the phenomena observed here needs to take judicious account of associated implicit selection effects.
We have excluded years in which no more than two hurricanes occur, as they supply no information about clustering within a year.
A \(k\)-means algorithm (\(k=20\), using Lloyd's algorithm based on \(10\) repetitions) classified the remaining 179 hurricanes from 37 years.
The 20 groups provided by the \(k\)-means algorithm are loosely ordered on an East-West axis as illustrated in Figure \ref{fig:WSKfig04} (see also Figure \ref{fig:WSKfiga4}, 
which uses a \(k\)-means analysis with \(k=4\), based on the entire 1851-2012 dataset).
In the Figure, the apparently anomalous barycentre trajectory running nearly horizontally near the more eastern end of the collection arises from a single rather long hurricane trajectory, 
whose initial behaviour 
 (including its up-crossing of \(20^\circ N\)) is removed as part of the process of confining
the barycentre trajectories to the strict intersection of registered time intervals for the component hurricane trajectories.
Consequently the grouping by \(k\)-means relates to trajectory behaviour on quite a narrow band of latitudes, as can be seen from Figure \ref{fig:WSKfig04}.
 \begin{table}
 \begin{small}
 \begin{tabular}{@{}cccccccccc@{}}
 \toprule
 Year & Labels & \hphantom{Labels}& \hphantom{Labels}& \hphantom{Labels}& \hphantom{Labels}& \hphantom{Labels}& \hphantom{Labels}& \hphantom{Labels}& \hphantom{Labels}\\
 \midrule
1950& 5& 1& 13& 8& 2& 13\\
1951& 12& 14& 7\\
1952& 6& 3& 8& 9\\
1953& 11& 12& 11& 4\\
1954& 11& 5& 5\\
1955& 7& 5& 13& 0& 14& 5& 11\\
1958& 5& 14& 7& 13& 8\\
1961& 14& 0& 17& 8& 8& 9\\
1962& 19& 14& 11\\
1963& 12& 13& 16& 12\\
1964& 7& 13& 2& 4& 12& 8& 3\\
1965& 1& 1& 17& 14\\
1966& 7& 9& 9\\
1969& 5& 13& 1\\
1975& 1& 13& 11\\
1976& 14& 19& 15\\
1979& 0& 11& 3& 1& 19& 16& 16\\
1980& 19& 19& 8& 17& 16\\
1981& 3& 18& 8& 13& 14& 13\\
1985& 1& 1& 5& 7\\
1988& 3& 0& 16& 10\\
1989& 12& 17& 17& 13& 5\\
1990& 14& 16& 15& 17\\
1995& 4& 8& 8& 16& 13& 11& 11& 19& 1\\
1996& 3& 8& 5& 11& 19& 15\\
1998& 7& 8& 17& 19& 19& 10\\
1999& 14& 3& 13& 3& 13\\
2000& 14& 4& 4& 14\\
2001& 11& 13& 19\\
2003& 1& 12& 5& 14\\
2004& 4& 6& 17& 2& 1& 2& 16& 17\\
2005& 1& 1& 1& 8& 15& 15& 10\\
2006& 3& 12& 14\\
2008& 12& 1& 0& 5& 0& 15\\
2010& 13& 5& 0& 12& 13\\
2011& 3& 8& 11& 13& 14\\
2012& 0& 13& 19& 12& 9\\
 \botrule
 \end{tabular}
 \end{small}
 \caption{179 hurricanes over 37 years, classified by year and by 20 groups using the \(k\)-means algorithm with \(k=20\).
Groups are ordered according to how westerly 
is the upcrossing by the corresponding barycentre trajectory of latitude \(35^\circ N\).\label{tab:labels}}{}
 \end{table}
Labelling the groups in this order, so that the most westerly group at first crossing of \(35^\circ N\) is given index 0, we obtain Table \ref{tab:labels}.
Reading from west to east, numbers in the \(20\) groups are given in Table \ref{tab:clusters}.
This procedure yields a test statistic (sum of numbers within each year of consecutive pairs belonging to the same group) of \(T=15\).
Computing mean and variance of \(T\), conditional on the numbers of each cluster occurring in each year, and assuming the normal approximation of \(T\) to be valid,
this can be referred to a conditional one-sided \(5\%\) level of \(10.66+4.15=14.81\). 
 \begin{table}
 \begin{small}
 \begin{tabular}{@{}cccccccccccccccccccc@{}}
 \toprule
  8   & 15   &  4   & 10   &  6   & 13   &  2   &  7   & 14   &  5   &  3   & 12   & 11   & 19   & 15   &  6   &  8   &  9   &  1   & 11   \\
 \botrule
\end{tabular}
 \end{small}
 \caption{Numbers of hurricanes in each of the 20 groups determined by the \(k\)-means algorithm with \(k=20\).
Groups are ordered according to how westerly is the upcrossing by the corresponding barycentre trajectory of latitude \(35^\circ N\).\label{tab:clusters}}{}
 \end{table}

\begin{figure}
 \includegraphics[width=10cm]{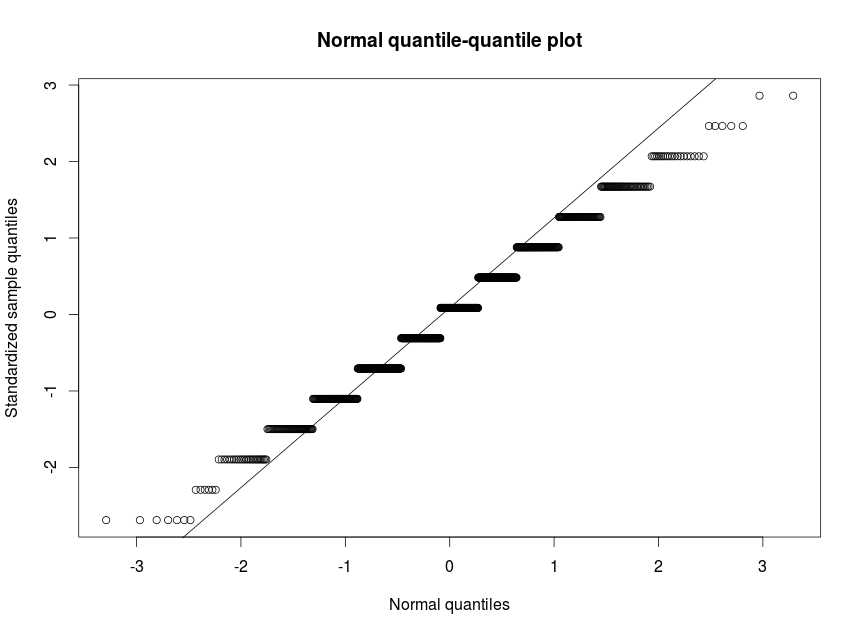}
 \caption{Quantile-quantile plot assessing approximate normality of the distribution of the test statistic \(T\)
 (with \(T\) constructed using \(k\)-means clustering with \(k=20\)) based on 1950-2012 hurricanes crossing latitudes \(20^\circ N\), \(35^\circ N\).}
 \label{fig:WSKfig05}
\end{figure}

A quantile-quantile plot, based on 1000 randomized versions of the data displayed in Table \ref{tab:labels}, 
shows that the standardized distribution of \(T\) has somewhat lighter tails when compared to a normal distribution (see Figure \ref{fig:WSKfig05}).
However a simulation test based on \(1000\) simulations yields an unremarkable \(p\)-value of \(7.5\%\), compared to a \(p\)-value of \(4.3\%\) using the normal approximation.

\begin{figure}
 \includegraphics[width=10cm]{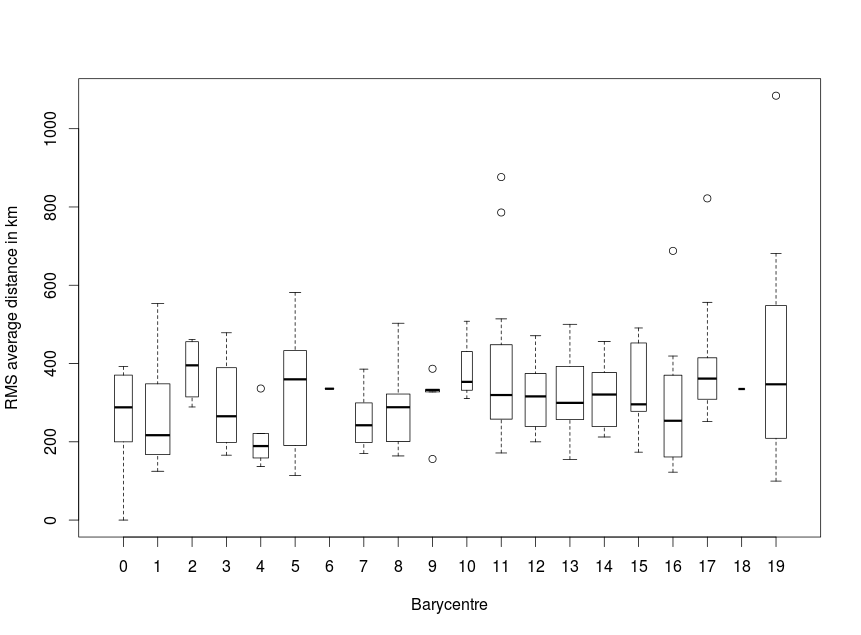}
 \caption{Boxplots of root-mean-square (RMS) average distances of hurricanes from associated barycentres (in kilometres).
 Boxplot widths are proportional to square-roots of sample sizes.}
 \label{fig:WSKfig08}
\end{figure}

If we replace the test statistic \(T\) by \(T_{\beta}\) (obtained by summing the contributions from \(T_{y,\beta}\) as defined in \eqref{eqn:potts-spread-stat}), 
so as to make a weighted count of repetitions at longer lags, then there is a slightly stronger indication of clustering. 
We choose \(\beta=0.25\), so that longer lags are penalized quite heavily, (results do not appear to be particularly sensitive to the choice of small \(\beta>0\).)
A quantile-quantile plot suggests normality of the distribution of  \(T_{\beta=0.25}\) under the hypothesis of no temporal association, though we omit this plot here,
and in any case focus on assessment \emph{via} a simulation-based permutation test.
We obtain \(T_{\beta=0.25}=16.97\), and a simulation test based on \(1000\) simulations yields a \(p\)-value of \(3.2\%\). 
(Examination of the data confirms that the dominant contribution to the modest improvement in \(p\)-value arises from pairs of similar hurricanes separated by lag \(2\).)
This therefore suggests modest evidence of temporal association among this particular set of hurricanes.
Given the hurricane trails have been classified into \(20\) clusters, over an East-West range of order \(6000\) km, the length scale of this association 
may be deemed to be of order of \(300\) km.
This is supported by the boxplots of root-mean-square average distances of hurricanes from associated barycentres in Figure \ref{fig:WSKfig08} 
(the calculation employs the approximation of Riemannian barycentres by cosine-barycentres). 
Mean distances between hurricanes and associated barycentres do indeed appear to be of order of 300 km. 
Visual inspection of the individual clusters confirms that some clusters do group together rather different hurricane trajectories,
underlining the need to be cautious in interpreting the formal statistical analysis given above.
We note the two groups containing fewer than 3 hurricanes (groups 6 and 18); inspection shows that trajectories in both groups exhibit atypical behaviour.
Various outliers and the larger dispersion of group 19
appear mostly to be linked to the more diverse behaviour of easterly hurricane trajectories.

\section{Conclusion}\label{sec:WSKconc}
This work illustrates how barycentres can be used in the analysis of trajectories with strong geometric content (here, hurricane trajectories lying on the surface of the terrestrial sphere).
The conclusions drawn are modest, namely that there is rather limited evidence in favour of temporal interaction. 
We repeat that potential selection effects need to be borne in mind here.
The nature of the dataset (secular trends, structure of temporal sequence of short time series) hamper further investigation.
Were the purpose of this paper to develop such an applied theme, then the next step would be to pay greater attention to other features of the underlying dataset (in particular, records of wind strength and atmospheric pressure),
and also to combine the above analysis with inference drawn from other associated meteorological datasets.
But the intention of this paper is more methodological: further development in such a direction could include the investigation of the more parametric inferential approaches mentioned in Section \ref{sec:WSKkmeans}:
including information derived from varying proportions of different groups of hurricanes, attempting maximum likelihood estimation of interaction parameters or even Bayesian inference exploiting the form of the likelihood of the heuristic parametric model. 

Other avenues of investigation would require commitment of more computational resource: investigation of the effect of the cosine-barycentre approximation, factoring out time-shifts, or making a discontinuous time-change by referring trajectories to their times of upcrossing of successive latitudes, or working in terms of arc-length rather than time.

Finally, and more speculatively, the geometric context of these data is very simple. Useful insight might arise from consideration of more ambitious questions. 
For example, and following one of the applications in \citet{SuDrydenKlassenLeSrivastava-2011}, this methodology could be extended to deal with the more complicated geometrical considerations that would arise 
when comparing three-dimensional trajectories arising in the study of chemotaxis, or more generally of motility of small organisms living at low Reynolds number \citep{Purcell-1977}.
One speculates that it might be possible to use these three-dimensional trajectories to draw inferences concerning stochastic characteristics of trajectories in the rotational group\index{rotational group}, 
using the imputed orientation of the small organism in question; at small length scales (where Brownian effects cannot be neglected) this might lead to intriguing statistical applications of 
the techniques underlying the celebrated Eells-Elworthy stochastic development\index{stochastic development} \citep[ch.VII.11]{Elworthy-1982}.


\bigskip

\noindent
Department of Statistics, University of Warwick, Coventry CV4 7AL, UK\\
Email: \url{w.s.kendall@warwick.ac.uk}\\
WWW: \url{www.warwick.ac.uk/wsk}

\bibliographystyle{chicago}

\bibliography{WSKMardia}

\end{document}